\documentclass[a4paper]{amsart}

\usepackage{amssymb,amscd,amsxtra}
\usepackage[all]{xy}

\theoremstyle{plain}
\newtheorem{thm}{Theorem}[section]
\newtheorem{lem}[thm]{Lemma}
\newtheorem{prop}[thm]{Proposition}
\newtheorem{cor}[thm]{Corollary}

\theoremstyle{definition}
\newtheorem{rem}[thm]{Remark}
\newtheorem{defn}[thm]{Definition}

\numberwithin{equation}{section}

\def\ie{\emph{i.e.}}
\def\ds{\displaystyle}
\def\:{\colon}
\def\.{\cdot}

\def\<{\left\langle}
\def\>{\right\rangle}
\def\({\left(}
\def\){\right)}
\def\ph#1{\phantom{#1}}
\def\epsilon{\varepsilon}
\def\subset{\subseteq}

\def\leq{\leqslant}
\def\geq{\geqslant}
\def\lla{\longleftarrow}

\def\lra{\longrightarrow}
\def\Lra{\Longrightarrow}

\def\tilde#1{\widetilde{#1}}
\def\iso{\cong}

\DeclareMathOperator{\rank}{rank}

\def\C{\mathbb{C}}

\def\F{\mathbb{F}}

\def\Q{\mathbb{Q}}
\def\N{\mathbb{N}}

\def\Z{\mathbb{Z}}
\def\B{\mathrm{B}}

\def\widebar#1{\overline{#1}}

\def\id{\mathrm{id}}

\def\CP{\C P}
\def\CPi{\CP^\infty}

\def\BP{\mathit{BP}}

\def\thh{\mathit{THH}}

\def\bp1{\BP\langle 1\rangle}
\def\qsym{{\sf QSymm}}
\def\sym{{\sf Symm}}
\def\nsym{{\sf NSymm}}
\def\lsc{\Omega \Sigma \mathbb{C}P^\infty}
\def\lse{\Omega \Sigma E}
\def\Lie{\mathrm{Lie}}

\def\clc{\cdot \ldots \cdot}
\DeclareMathOperator{\Prim}{Prim}

\begin{document}
\title[$\qsym$ -- a topological point of view]
{Quasisymmetric functions from a topological point of view}
\author{Andrew Baker \and Birgit Richter}
\address{Department of Mathematics, University of Glasgow,
Glasgow G12 8QW, Scotland.}
\email{a.baker@maths.gla.ac.uk}
\urladdr{http://www.maths.gla.ac.uk/$\sim$ajb}
\address{Department Mathematik der Universit\"at Hamburg,
20146 Hamburg, Germany.}
\email{richter@math.uni-hamburg.de}
\urladdr{http://www.math.uni-hamburg.de/home/richter/}
\thanks{The authors thank the Mittag-Leffler institute for support and
hospitality allowing major parts of this work to be carried out during
their stay there. The first-named author was supported by the Carnegie
Trust for the Universities of Scotland.}

\begin{abstract}
It is well-known that the homology of the classifying space of the unitary
group is isomorphic to the ring of symmetric functions $\sym$. We offer
the cohomology of the space $\lsc$ as a topological model for the ring
of quasisymmetric functions $\qsym$. We exploit standard results from
topology to shed light on some of the algebraic properties of $\qsym$. In
particular, we reprove the Ditters conjecture. We investigate a product
on $\lsc$ that gives rise to an algebraic structure which generalizes the
Witt vector structure in the cohomology of $BU$. The canonical Thom spectrum
over $\lsc$ is highly non-commutative and we study some of its features,
including the homology of its topological Hochschild homology spectrum.
\end{abstract}

\maketitle

\section*{Introduction}

Let us recall some background on the variants of symmetric functions.
For a much more detailed account on that see~\cite{H1,H2}.

The algebra of symmetric functions, $\sym$, is an integral graded
polynomial algebra
\[
\sym = \Z[c_1,c_2,\ldots],
\]
where $c_i$ has degree $2i$. The reader is encouraged to think of
these $c_i$ as Chern classes. This algebra structure can be extended
to a Hopf algebra structure by defining the coproduct to be that
given by the Cartan formula
\[
\Delta(c_n) = \sum_{p+q=n} c_p \otimes c_q.
\]
The antipode on $\sym$ is defined as
\[
\chi(c_n) = \sum_{i_1+\ldots+i_m=n} (-1)^{m} c_{i_1}\clc c_{i_m}.
\]
This Hopf algebra is self dual in the sense that there is an isomorphism
of Hopf algebras
\[
\sym^* \iso \sym,
\]
where $\sym^*$ is the degree-wise $\Z$-linear dual of $\sym$. In
particular, $\sym$ is bipolynomial, \ie, the underlying algebra
of the Hopf algebra and its dual are both polynomial algebras.

The non-commutative analogue of the algebra $\sym$ is the algebra
of non-symmetric functions, $\nsym$ (also known as the Leibniz
algebra) which is the free associative graded $\Z$-algebra
$\Z\< Z_1, Z_2, \ldots \>$ on generators $Z_1,Z_2,\ldots$, where
$Z_i$ has degree~$2i$. Again, $\nsym$ comes with a natural
coproduct given by
\[
\Delta(Z_n) = \sum_{p+q=n} Z_p \otimes Z_q
\]
and an antipode given on $Z_n$ by
\[
\chi(Z_n) =
\sum_{a_1+\ldots+a_m = n} (-1)^{m}Z_{a_1}\cdot \ldots \cdot Z_{a_m}.
\]

The Hopf algebra of \emph{quasisymmetric functions} (sometimes
written \emph{quasi\--sym\-met\-ric functions}), $\qsym$, is defined
to be the dual Hopf algebra to $\nsym$. We follow the convention
from~\cite{H1,H2}, denoting the element dual to the monomial
$Z_{a_1} \clc Z_{a_n}$ with respect to the monomial basis by
$\alpha = [a_1,\ldots, a_n]$ and call the number $a_1 + \ldots + a_n$
the \emph{degree} of $\alpha$. The resulting product structure among
these elements is given by the \emph{overlapping shuffle product}
of~\cite[section~3]{H01}. For example, using the dual pairing, we
find that
\[
[3][1,2] = [3,1,2] + [1,3,2] + [1,2,3]+ [4,2] + [1,5].
\]

Often it is useful to vary the ground ring and replace the integers by
some other commutative ring with unit $R$. We define $\sym(R)$ to be
$\sym \otimes R$, and similarly we set $\nsym(R) = \nsym \otimes R$
and $\qsym(R) = \qsym \otimes R$.

The algebras $\nsym$ and $\qsym$ have received a great deal of attention
in combinatorics. Several structural properties were proven, for instance
about the explicit form of the primitives in the coalgebra $\nsym$~\cite{H3}
or the freeness of $\qsym$ as a commutative algebra~\cite{H01}. The latter
result is known as the \emph{Ditters conjecture}, and is our
Theorem~\ref{thm:QSymm-Poly}. The original methods of proof came from
within combinatorics. We offer an alternative proof using ingredients from
algebraic topology.

In the case of symmetric functions, Liulevicius~\cite{L} exploited the
identification of $\sym$ with the cohomology of $BU$ to use topology to
aid the understanding of some of the properties of $\sym$. In this paper
we offer a topological model for the Hopf algebras $\nsym$ and $\qsym$
by interpreting them as homology and cohomology of the loop space of
the suspension of the infinite complex projective space, $\lsc$.

Our desire to find a topological model for the algebra of quasisymmetric
functions has its origin in trying to understand Jack Morava's thoughts 
on connecting Galois theory of structured ring spectra to motivic Galois
theory, as explained in~\cite{Mo}. We do not claim that our insights are 
helpful in this context, but it was our motivation to start this 
investigation.

\bigskip

The first part of this paper is concerned with the algebraic structure of 
the algebra of quasisymmetric functions.

In Section~\ref{sec:HAstructure} we describe the isomorphism between the
Hopf algebra $\nsym$ and the homology of the loop space on the suspension
of the infinite complex projective space. This identification is probably
known to many people, but we do not know of any source where this is
seriously exploited.

We give a proof of the Ditters conjecture in Section~\ref{sec:qsympoly}.
This conjecture states that the algebra $\qsym$ is polynomial and was
proven by Hazewinkel in \cite{H01}. Our proof uses the Hilton-Milnor
theorem which also yields an explicit set of generators over the rationals.

\medskip

The second part deals with the topological properties of the model
$\Omega \Sigma \CPi$ and its relation to $BU$.

We investigate the $p$-local structure of $\qsym$ in Section~\ref{sec:plocal}
using the splitting of $\Sigma \CPi$ at a prime~$p$. We discuss Steenrod
operations on $\qsym(\F_p)$ in Section~\ref{sec:steenrod}.

It is well-known that the ring of big Witt vectors on a commutative
ring is represented by $\sym$. Topologically this structure is induced
by the two canonical $H$-structures on $BU$. We recall this in
Section~\ref{sec:Witt&LambdaRings}, then in Section~\ref{sec:q-Witt}
we introduce product structures on $\lsc$ which in cohomology produces
a structure that we call the \emph{quasi-Witt} structure  on $\qsym$,
a hitherto unremarked algebraic structure, which differs from the one
explored by Hazewinkel in~\cite[\S 14]{H1}.

The canonical map from $\lsc$ to $BU$ is a loop map. Therefore the
associated Thom spectrum has a strictly associative multiplication.
But as is visible from the non-commutativity of its homology, it is
not even homotopy commutative. We will describe some of its features
in Sections~\ref{sec:ThomSpectrum-xi} and investigate the homology of
its topological Hochschild homology spectrum in Section~\ref{sec:thh}.

\part{Algebraic properties of the algebra of quasisymmetric functions}

\section{A topological manifestation of $\nsym$}
\label{sec:HAstructure}

In this part we will recall some standards facts about $H_*(\lsc)$.

There is a nice combinatorial model $JX$ for any topological space of the
form $\Omega \Sigma X$ with $X$ connected, namely the James construction
on $X$, fully described in~\cite[VII \S 2]{Wh}. After one suspension
this gives rise to a splitting
\begin{equation}\label{eqn:splitting}
\Sigma \Omega \Sigma X \sim \Sigma JX \sim \bigvee_{n\geq1}\Sigma X^{(n)},
\end{equation}
where $X^{(n)}$ denotes the $n$-fold smash power of $X$.

If the homology of $X$ is torsion-free, then the homology of $JX$
is the tensor algebra on the reduced homology of $X$
\[
H_*(JX) \cong T(\tilde{H}_*(X)).
\]
The concatenation of loops in $\lsc$ together with the diagonal
on $\lsc$ turns the homology of $\lsc$ into a Hopf algebra.

The integral homology of $\mathbb{C}P^\infty$ has
$H_i(\mathbb{C}P^\infty) = \mathbb{Z}$ for all even $i$ with generators
$\beta_i \in H_{2i}(\CPi)$ and is trivial in odd degrees. Therefore
\[
H_*(\lsc) \cong T(\tilde{H}_*(\mathbb{C}P^\infty))
                      = \Z\langle \beta_1,\beta_2,\ldots \rangle,
\]
with $\beta_i$ being a non-commuting variable in degree $2i$. Thus
there is an isomorphism of algebras
\begin{equation}\label{eqn:homology}
H_*(\lsc) \cong \nsym
\end{equation}
under which $\beta_n$ corresponds to $Z_n$. The coproduct $\Delta$
on $H_*(\lsc)$ induced by the diagonal in $\lsc$ is compatible with
the one on $\nsym$:
\[
\Delta(\beta_n) = \sum_{p+q=n} \beta_p \otimes \beta_q.
\]
Putting this information together, we see that~\eqref{eqn:homology}
gives an isomorphism of graded, connected Hopf algebras. Note that
the antipode $\chi$ in $H_*(\lsc)$ arises from the time-inversion
of loops. As antipodes are unique for Hopf algebras which are
commutative or cocommutative, this gives a geometric interpretation
for the antipode in $\nsym$.

As the homology of $\lsc$ is a graded free abelian group, the linear
dual of $H_*(\lsc)$ is canonically isomorphic to the cohomology, $H^*(\lsc)$,
which is also a Hopf algebra. Thus we have proven the following result.
\begin{thm}\label{thm:qsym}
There are isomorphisms of graded Hopf algebras
\[
H_*(\lsc) \cong \nsym,
\quad
H^*(\lsc) \cong \qsym.
\]
\end{thm}
\begin{rem}\label{rem:degrees}
Note that the cohomological degree of a generator corresponding
to a sequence $\alpha = [a_1,\ldots,a_n]$ is twice its degree.
\end{rem}

There is the canonical inclusion map $j\:\CPi = BU(1) \lra BU$.
The universal property of $\lsc$ as a free $H$-space gives an
extension to a loop map
\[
\xymatrix{
{\lsc} \ar@{.>}[dr]^{j} & {} \\
{\CPi} \ar[r] \ar[u]& {BU}
}
\]
which induces a virtual bundle $\xi$ on $\lsc$. On homology, the
map~$j$ induces an epimorphism because $\CPi$ gives rise to the
algebra generators in $H_*(BU)$ and correspondingly, $H^*(j)$ is
a monomorphism on cohomology. This corresponds to the inclusion
of symmetric functions into quasisymmetric functions. We will
describe this in more detail in Section~\ref{sec:Witt&LambdaRings}.

Hazewinkel mentions in~\cite{H3} that over the rationals the Lie
algebra of primitives in the Hopf algebra $\nsym(\Q)$ is free and
says that the primitive part of $\nsym$ ``is most definitely not
a free Lie algebra; rather it tries to be something like a divided
power Lie algebra (though I do not know what such a thing would be)''.
In this section we give a topological proof of the rational result
and we explain how to make sense of this last comment in positive
characteristic.

A theorem of Milnor and Moore~\cite[Appendix]{M&M} identifies the
Lie-algebra of primitives in the Hopf algebra $H_*(\lsc; \mathbb{Q})$
with the Lie-algebra $\pi_*(\lsc) \otimes \Q$. Here the Lie-algebra
structure on $\pi_*(\lsc)$ is that given by the Samelson-Whitehead
product~\cite[X \S\S 5-7]{Wh}.

Let $C$ be a simply connected rational co-$H$-space that is a CW-complex.
Scheerer in~\cite[pp.72--73]{Sch} proves that for such spaces $C$, the
Lie algebra $\pi_*(\Omega C)$ is a free Lie-algebra. So in particular,
$\pi_*(\lsc) \otimes \Q$ is a free Lie-algebra.

In the case of positive characteristic, operads help to identify the
primitives in $\nsym$. In~\cite[theorem~1.2.5]{F}, Fresse uses the
fact that for a vector space $V$ over a field $k$ of characteristic~$p$,
the primitives in a tensor algebra $T(V)$ can be identified with the
free $p$-restricted Lie-algebra generated by $V$. He shows that the
free $p$-restricted Lie-algebra is isomorphic to the direct sum of
invariants
\[
\bigoplus_{n\geq 1} (\Lie(n) \otimes V^{\otimes n})^{\Sigma_n},
\]
where $\Lie(n)$ is the $n$th part of the operad which codifies Lie
algebras. We note that instead of invariants we may take coinvariants
\[
\bigoplus_{n\geq 1} (\Lie(n) \otimes V^{\otimes n})_{\Sigma_n}
\]
to give the free Lie algebra generated by the vector space $V$, whereas
in Fresse's terminology of~\cite[p.4122]{F}, the invariants codify the
\emph{free Lie algebra with divided symmetries}. In the case when~$V$
is the vector space generated by $Z_1,Z_2,\ldots$ we may deduce the
following result.
\begin{prop}\label{prop:FreeLie}
For a field $k$ of positive characteristic $p$, the Lie subalgebra
of primitives $\Prim(\nsym(k))$ agrees with the free $p$-restricted
Lie-algebra on the $k$-vector space $V$ generated by $Z_1,Z_2, \ldots $
and furthermore there is an isomorphism
\[
\bigoplus_{n\geq 1} (\Lie(n) \otimes V^{\otimes n})^{\Sigma_n}
\iso
\Prim(\nsym(k)).
\]
\end{prop}

\section{A proof of the Ditters conjecture}
\label{sec:qsympoly}

In this section we give a topological proof of the Ditters conjecture
which asserts that the algebra $\qsym$ is a free commutative algebra
(see~\cite{H01, H2}), and use the Hilton-Milnor theorem to show that
over the rationals the generators can be indexed on Lyndon words (see
Definition~\ref{defn:Lyndon}). The Ditters conjecture started off as
a statement~\cite[proposition~ 2.2]{D}, but it turned out that the proof
was not correct. There were later attempts by Ditters and Scholtens to
prove the conjecture~\cite{DS}, however, the line of argument there
turned out to be incorrect as well. Hazewinkel proved the conjecture
in~\cite{H01}. For another approach on related matters from a topological
perspective see~\cite{C}. Here is our statement of these results.
\begin{thm}\label{thm:QSymm-Poly}
The algebra of quasisymmetric functions, $\qsym$, is a free commutative
algebra. Over the rationals, the polynomial generators of $\qsym(\Q)$
in degree~$2n$  can be indexed on Lyndon words of degree~$n$.
\end{thm}

This Theorem recovers Hazewinkel's result~\cite[theorem 8.1]{H01}.
Recall that the degree of a word $a_1 \clc a_n$ with $a_i \in \N$
is $a_1+ \ldots + a_n$.

Our proof proceeds by using Borel's theorem on the structure of Hopf
algebras over perfect fields~\cite[theorem 7.11]{M&M} to first identify
the rationalization of $\qsym$ as a polynomial algebra and then to show
that the $\F_p$-reductions are polynomial for all primes~$p$. Finally,
we use a gluing result Proposition~\ref{prop:PolyAlgs} to obtain the
integral statement. For the explicit form of the generators we compare
the words arising in the Hilton-Milnor theorem to Lyndon words in
Proposition~\ref{prop:Algorithm}.

The rational version of Borel's theorem immediately implies that the
algebra of rational quasisymmetric functions, $\qsym(\Q)$, is a polynomial
algebra because all its generators live in even degrees.

Rationally the suspension of $\CPi$ splits into a wedge of rational
spheres
\[
\Sigma \CPi_\Q \sim \;\mathbb{S}^3_\Q \vee \mathbb{S}^5_\Q \vee \ldots
 \; \sim \;\Sigma(\mathbb{S}^2_\Q \vee \mathbb{S}^4_\Q \vee \ldots).
\]
Therefore we can apply the Hilton-Milnor theorem for loops on the
suspension of a wedge of spaces~\cite[theorem~1.2]{T} and obtain
\[
\Omega \Sigma(\mathbb{S}^2_\Q \vee \mathbb{S}^4_\Q \vee \ldots)
                       \sim \prod_{\alpha}\!{}^{'} \Omega N_\alpha,
\]
where after suitable suspension, $N_\alpha$ is a smash product of
rational spheres, thus its cohomology is monogenic polynomial. The
$\alpha$ in the indexing set of the weak product run over all
\emph{basic products} in the sense of Whitehead~\cite[pp.~511--512]{Wh}.

Note that in the usual formulation of the Hilton-Milnor theorem, only
a finite number of wedge summands are considered. However, a colimit
argument gives the countable case as well. In
Subsection~\ref{subsec:lyndon} we give a bijection between the set of
basic products and the set of Lyndon words. Thus
\[
\qsym(\Q) \cong \Q[x_\alpha\;|\; \text{$\alpha$ Lyndon}].
\]

Now we consider the $\F_p$-cohomology $H^*(\lsc;\F_p)$ which is canonically
isomorphic to the mod~$p$ reduction of the
quasisymmetric functions, $\qsym(\F_p)$. A priori the Borel theorem allows
for truncated polynomial algebras. However we will use the action of the
Steenrod algebra to prove:
\begin{prop}\label{prop:QSymm-modp}
$\qsym(\F_p)$ is a polynomial algebra.
\end{prop}
\begin{proof}
Using the James splitting we obtain that as a module over the Steenrod
algebra $\mathcal{A}_p^*$, the positive part of $\qsym(\F_p)$ has a
direct sum decomposition
\[
\qsym(\F_p)^{*>0} \cong \bigoplus_n \tilde{H}^*((\CPi)^{(n)};\F_p).
\]
We have to show that no $p$th power of an element $x$ can be zero. Such
a power corresponds to the Steenrod operation  $\mathcal{P}^{|x|/2}$
applied to $x$. If $p$ is odd, we write $\mathcal{P}^i$ for the reduced
power operation, while for $p=2$, we set $\mathcal{P}^i=Sq^{2i}$.

These operations are non-trivial on the cohomology of $\CPi$ and from
the Cartan formula we see that they are non-trivial on the cohomology
of the smash powers.
\end{proof}

\begin{rem}\label{rem:QSymm-modp}
The above proof showed that the $p$th power operation on the algebra
$\qsym(\F_p)$ is non-trivial. We will explicitly determine the action
of the mod~$p$ Steenrod algebra on $\qsym(\F_p)$ in Section~\ref{sec:steenrod}.
\end{rem}

Taking the rational result and the $\F_p$-cohomology result together
with Proposition~\ref{prop:PolyAlgs} below yields a proof of
Theorem~\ref{thm:QSymm-Poly}.

\subsection{An auxiliary result on polynomial algebras}
\label{sec:PolyAlgs}

In this section we provide a useful local to global result on polynomial
algebras which may be known but we were unable to locate a specific
reference.

If $R$ is a commutative ring, then for a graded $R$-algebra  $A^*$, we
write $DA^n$ (resp. $QA^n$) for the decomposables (resp. indecomposables)
in degree~$n$. Unspecified tensor products will be taken over whatever
ground ring $R$ is in evidence. If $p$ is a positive prime or $0$, let
$\F_p$ denote either the corresponding Galois field or $\F_0=\Q$.

\begin{prop}\label{prop:PolyAlgs}
Let $H^*$ be a graded commutative connective $\Z$-algebra which is
concentrated in even degrees and with each $H^{2n}$ a finitely generated
free abelian group. If for each non-negative rational prime~$p$,
$H(p)^*=H^*\otimes\F_p$ is a polynomial algebra, then $H^*$ is a polynomial
algebra and for every non-negative rational prime~$p$,
\[
\rank QH^{2n} = \dim_{\F_p} QH(p)^{2n}.
\]
\end{prop}
\begin{proof}
Let $p\geq0$ be a prime. We will denote by $\pi_p(n)$ the number of polynomial
generators of $H(p)^*$ in degree $2n$, and this is equal to $\dim_{\F_p} QH(p)^{2n}$.
The Poincar\'e series of the even degree part of $H^*(p)$ satisfies
\[
\sum_{n\geq0}\rank H^{2n} t^n = \sum_{n\geq0}\rank H(p)^{2n} t^n
                              = \prod_{n\geq0}(1-t^n)^{-\pi_p(n)},
\]
hence $\pi_p(n)$ is independent of~$p$.

Now it is easy to see that the natural homomorphism
$DH^{2n}\otimes\Q \lra DH(0)^{2n}$ is an isomorphism and therefore
the natural homomorphism $QH^{2n}\otimes\Q \lra QH(0)^{2n}$ is an
isomorphism. Furthermore, for each positive prime~$p$ there is a
commutative diagram
\[
\xymatrix{
{0} \ar[r] & {DH^{2n}} \ar[r] \ar[d]^{\mathrm{epic}}& {H^{2n}} \ar[r]
\ar[d]^{\mathrm{epic}} &  {QH^{2n}} \ar[r] \ar[d] & 0 \\
{0} \ar[r] & {DH(p)^{2n}} \ar[r] & {H(p)^{2n}} \ar[r] & {QH(p)^{2n}}
\ar[r] & {0}
}
\]
with exact rows and whose columns have the indicated properties. Since
the right hand vertical homomorphism factors as
\[
QH^{2n} \lra QH^{2n}\otimes \F_p \lra QH(p)^{2n},
\]
we see that the right hand factor is an epimorphism
$QH^{2n}\otimes \F_p \lra QH(p)^{2n}$. This implies that
\[
\pi_0(n) \geq \pi_p(n).
\]

We will now show that the indecomposables in degree $2n$ are torsion
free. Assume that $QH^{2n}$ were of the form
\[
QH^{2n} = \mathbb{Z}^r \oplus \bigoplus_{p} T_p,
\]
where $p$ runs over a finite set of primes and $T_p$ is the $p$-torsion
subgroup of $QH^{2n}$. Then
\[
QH^{2n} \otimes \F_p \cong \F_p^r \oplus T_p \otimes \F_p.
\]
Thus if $T_p \neq 0$, then the dimension of  $QH^{2n} \otimes \F_p$ as
an $\F_p$-vector space would be strictly bigger than~$r$. However, the
following argument shows that these dimensions are equal and so $T_p$
has to be zero.

Let $I = H^{* > 0}$ and $I_{\F_p} = H^{*>0} \otimes \F_p = I \otimes \F_p$.
Observe that
\[
I \otimes \F_p \otimes I \otimes \F_p \cong (I \otimes I) \otimes \F_p
\]
and consider the following commuting diagram with exact rows.
\[
\xymatrix{
{(I \otimes I) \otimes \F_p} \ar[rr]^{\ph{-}\mathrm{mult} \otimes \id}
                                                     \ar[d]^{\cong} & &
{I \otimes \F_p} \ar[r] \ar[d]^{\cong} & {I/I^2 \otimes \F_p}
                                            \ar[r] \ar@{.>}[d] & {0} \\
{I_{\F_p} \otimes I_{\F_p}} \ar[rr]^{\ph{-}\mathrm{mult}}
            & & { I_{\F_p}} \ar[r]& { I_{\F_p}/I_{\F_p}^2} \ar[r] & {0}
}
\]
Thus we obtain that $I/I^2 \otimes \F_p \cong I_{\F_p}/I_{\F_p}^2$,
\ie, $QH^{2n}\otimes \F_p \cong  QH(p)^{2n}$.

Now for each $n$, choose a lifting of a basis of $QH(p)^{2n}$ to linearly
independent elements $x_{n,i}$ of $H^{2n}$. It is clear that under the
natural monomorphism $H^{2n} \lra H(0)^{2n}$, these give a part of a
polynomial generating set for $H(0)^{2n}$ and therefore generate a
polynomial subalgebra $P^* = \Z[x_{n,i}:n,i] \subset H^*$. We will use
induction on degree to show that we have equality here.

For $n=1$, we have $QH^2=H^2$. Now suppose that we have $H^* = P^*$ in
degrees less than $2k$. Then $DH^{2k} = DP^{2k}$ and for each positive
prime~$p$ there is a diagram with short exact rows
\[
\xymatrix{
{0} \ar[r] & {DP^{2k}} \ar[r] \ar[d]^{=} & {P^{2k}} \ar[r]
\ar[d]^{\mathrm{incl}}& {QP^{2k}} \ar[r] \ar[d] & {0}  \\
{0} \ar[r]  & {DH^{2k}}  \ar[r]
\ar[d]^{\mathrm{epic}} & {H^{2k}} \ar[r] \ar[d]\ar[d]^{\mathrm{epic}}
& {QH^{2k}}  \ar[r]
\ar[d]^{\mathrm{epic}}& {0} \\
{0} \ar[r] & {DH(p)^{2k}} \ar[r] & {H(p)^{2k}} \ar[r] &
{QH(p)^{2k}} \ar[r] & {0}
}
\]
in which the composite homomorphism $P^{2k} \lra QH(p)^{2k}$ is
surjective. To complete the proof, we must show that the cokernel
of the inclusion $P^{2k} \lra H^{2k}$ is trivial. Since $P^{2k}$
and $H^{2k}$ agree rationally, this cokernel is torsion and it
suffices to verify this locally at each prime~$p$. The map from
$P^{2k} \otimes \F_p$ to $H(p)^{2k}$ is an isomorphism, so we
see that over the local ring $\Z_{(p)}$,
\[
P^{2k}_{(p)}+pH^{2k}_{(p)} = H^{2k}_{(p)}
\]
and Nakayama's Lemma implies that $P^{2k}_{(p)} = H^{2k}_{(p)}$.
Thus $P^{2k} = H^{2k}$ and we have established the induction step.
\end{proof}

\subsection{Basic products and Lyndon words}
\label{subsec:lyndon}

Usually~\cite{D, H01} the set of polynomial generators of the rationalized
algebra of quasisymmetric functions is indexed on Lyndon words, whereas
our approach yields a polynomial basis indexed on basic products in the
sense of~\cite[pp.~511--512]{Wh}. The aim of this section is to compare
these two sets of generators.

First, let us recall some notation and definitions. See~\cite[\S 5]{Re}
for more details.

Let $A$ be an alphabet, finite or infinite. We assume that $A$ is linearly
ordered. The elements of $A$ are called \emph{letters}. Finite sequences
$a_1\ldots a_n$ with $a_i \in A$ are \emph{words}; the number of letters
in a word is its length. We use the lexicographical ordering on words given
as follows. A word $u$ is smaller than a word $v$ if and only if $v = ur$,
where $r$ is a non-empty word or if $u= wau'$ and $v = wbv'$ where $w,u',v'$
are words and $a$ and $b$ are letters with $a < b$. If a word $w$ can be
decomposed as $w = uv$, then $v$ is called a right and $u$ is called the
left factor. If $u$ is not the empty word, then $v$ is called a proper
right factor.

\begin{defn}\label{defn:Lyndon}
A word is \emph{Lyndon} if it is a non-empty word which is smaller than
any of its proper right factors.
\end{defn}

So for example, if $A$ is the alphabet that consists of the natural numbers
with its standard ordering, then the first few Lyndon words are
\[
1,2,3,\ldots, 12, 13, 23,\ldots,112, \ldots.
\]

Words are elements in the free monoid generated by the alphabet $A$. Let
$M(A)$ be the free magma generated by $A$, \ie, we consider non-associative
words built from the letters in $A$. Elements in the free magma correspond
to planar binary trees with a root where the leaves are labelled by the
elements in $A$. For instance, if $a,b,c$ are letters, then the element
$(ab)c$ corresponds to the tree

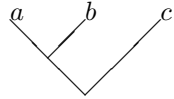
\vspace{1cm}
\begin{center}
\begin{picture}(2,3)
\setlength{\unitlength}{1cm}
\put(1,0){\line(-1,1){1}}
\put(1,0){\line(1,1){1}}
\put(0.5,0.5){\line(1,1){0.5}}
\put(0,1){$a$}
\put(1,1){$b$}
\put(2,1){$c$}
\end{picture}
\end{center}
whereas $a(bc)$ corresponds to

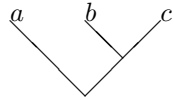
\vspace{1cm}
\begin{center}
\begin{picture}(2,3)
\setlength{\unitlength}{1cm}
\put(1,0){\line(-1,1){1}}
\put(1,0){\line(1,1){1}}
\put(1.5,0.5){\line(-1,1){0.5}}
\put(0,1){$a$}
\put(1,1){$b$}
\put(2,1){$c$}
\end{picture}
\end{center}

We will briefly recall the construction of basic products. For more
details see~\cite{Wh}.

Let $A$ be an alphabet whose letters are linearly ordered, for instance
$A = \{a_1,\ldots, a_k\}$ with $a_1 < \ldots < a_k$. Basic products of
length one are just the letters $a_i$. We assign a \emph{rank} and a
\emph{serial number} to each basic product. The convention for basic
products of length one is that the serial number of $a_i$ is $s(a_i) = i$,
whereas the rank is  $r(a_i) = 0$.

Assume that basic products of length up to $n-1$ have been already defined
and that these words are linearly ordered in such a way that a word $w_1$
is less than a word $w_2$ if the length of $w_1$ is less than the length
of $w_2$, and assume that we have assigned ranks to all those words. Then
the basic products of length $n$ are all (non-associative) words of length~$n$
of the form $w_1w_2$ such that the $w_i$ are basic products, $w_2 < w_1$
and the rank of $w_1$, $r(w_1)$, is smaller than the serial number of $w_2$,
$s(w_2)$. We choose an arbitrary linear ordering on these products of length~$n$
and we define their rank as $r(w_1w_2)= s(w_2)$.

Any basic product on $A$ is in particular an element of the free magma
generated by $A$. We need the fact that on an alphabet with $k$ elements,
the number of basic products of length~$n$ is
\[
\frac{1}{n}\sum_{d\mid n}\mu(d)k^{n/d},
\]
where $\mu$ is the M\"obius function~\cite[p.~514]{Wh}. This number
agrees with the number of Lyndon words of length $n$ on such an
alphabet, see~\cite[theorem~5.1, \& corollary 4.14]{Re} for details.
So we know that there is an abstract bijection between the two sets.

\begin{prop}\label{prop:Algorithm}
There is a canonical algorithm defining a bijection between the set
of basic products and the set of Lyndon words on a finite alphabet~$A$.
\end{prop}

Roughly speaking, the idea of the proof is to `correct' the ordering
in the product $w_1w_2$ of basic products with $w_2 < w_1$ and to switch
the product back to one that is ordered according to the convention used
for building Lyndon words.

\begin{proof}
We start with a basic product $\xi$ of length $n$ and consider it as
an element in the magma generated by $A$ and take its associated planar
binary tree. A binary tree has a natural level structure: we regard a
binary tree such as

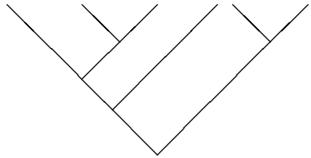
\vspace{2cm}\hspace{3cm}
\begin{picture}(4,2)
\setlength{\unitlength}{1cm}
\put(2,0){\line(-1,1){2}}
\put(2,0){\line(1,1){2}}
\put(1.4,0.6){\line(1,1){1.4}}
\put(1,1){\line(1,1){1}}
\put(1.5,1.5){\line(-1,1){0.5}}
\put(3.5,1.5){\line(-1,1){0.5}}
\end{picture}

\par \medskip \noindent
as having four levels

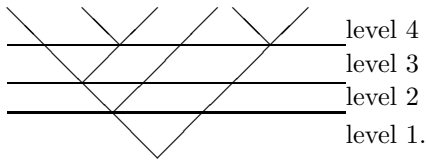
\vspace{2cm} \hspace{3cm}
\begin{picture}(5,2)
\setlength{\unitlength}{1cm}
\put(2,0){\line(-1,1){2}}
\put(2,0){\line(1,1){2}}
\put(1.4,0.6){\line(1,1){1.4}}
\put(1,1){\line(1,1){1}}
\put(1.5,1.5){\line(-1,1){0.5}}
\put(3.5,1.5){\line(-1,1){0.5}}
\put(0,1.5){\line(1,0){4.5}}
\put(0,1){\line(1,0){4.5}}
\put(0,0.6){\line(1,0){4.5}}
\put(4.5,1.6){level 4}
\put(4.5,1.15){level 3}
\put(4.5,0.7){level 2}
\put(4.5,0.2){level 1.}
\end{picture}

Starting with a basic product of length $n$, its tree has some
number of levels, say $m$. The idea is to work from the leaves
of the tree to its root and transform the basic product into
a Lyndon word during this process.

We start with level $m$. Whenever there is a part of a tree of
{\sf V}-shape we induce a multiplication; otherwise we leave the
element as it is. In the above example, the starting point could
be a word such as $((a_1(a_2a_3))a_4)(a_5a_6)$.

{}From level $m$ to level $m-1$ we induce multiplication whenever
there is a local picture like

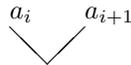
\vspace{1cm}
\begin{center}
\begin{picture}(2,2)
\setlength{\unitlength}{1cm}
\put(1,0){\line(-1,1){0.5}}
\put(1,0){\line(1,1){0.5}}
\put(0.5,0.6){$a_i$}
\put(1.5,0.6){$a_{i+1}$}
\end{picture}
\end{center}
and in such a case we send $(a_ia_{i+1})$ to
\[
\begin{cases}
a_ia_{i+1} & \text{if $a_i < a_{i+1}$},  \\
a_{i+1}a_{i} & \text{if $a_{i} > a_{i+1}$}.
\end{cases}
\]
Note that in a basic product equal neighbours $a_i = a_{i+1}$ do
not occur. We place this product on the corresponding leaf in
level $m-1$.

Iterating this procedure gives a reordered word in the alphabet~$A$.

We have to prove that this is a Lyndon word. We do this by showing
that in each step in the algorithm we produce Lyndon words. Starting
with the highest level $m$ this is clear because words of length two
of the form $ab$ with $a < b$ are Lyndon words.

In the following we will slightly abuse notation and use $\prod_{j=1}^t a_i$
for the ordered product $a_1\cdot \ldots \cdot a_t$. Assume that we have
reduced the tree down to an intermediate level less than $m$ and obtained
Lyndon words as labels on the leaves. In the next step we have to check
that every multiplication on subtrees of the form

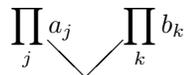
\vspace{1cm}
\begin{center}
\begin{picture}(2,2)
\setlength{\unitlength}{1cm}
\put(1,0){\line(-1,1){0.5}}
\put(1,0){\line(1,1){0.5}}
\put(0,0.6){$\ds\prod_j a_j$}
\put(1.5,0.6){$\ds\prod_k b_k$}
\end{picture}
\end{center}
with $\prod_{j=1}^t a_j$ and $\prod_{k=1}^s b_k$ Lyndon words, again
give a Lyndon word. But this is proved in~\cite[(5.1.2), p.106]{Re}.

Each  step in the algorithm is reversable, therefore the algorithm
defines an injective map. As the cardinalities of the domain and
target agree, this map is a bijection.
\end{proof}

\part{$\mathbf{\Omega \Sigma \CPi}$: a splitting, Witt vectors and
its associated Thom spectrum}

\section{A $p$-local splitting}
\label{sec:plocal}

In this section we fix an odd prime $p$. On the one hand,
from~\cite[lecture 4]{Ad-0}
we have the Adams splitting of $BU$ localized at~$p$,
\begin{equation}\label{eqn:Asplitting}
BU_{(p)} \simeq W_1 \times \ldots \times W_{p-1}
\end{equation}
which is a splitting of infinite loop spaces. On the other hand,
by~\cite{McG} there is an unstable $p$-local splitting of
$\Sigma\CPi_{(p)}$ into a wedge of spaces,
\begin{equation}\label{eqn:Csplitting}
\Sigma \CPi_{(p)} \simeq Y_1 \vee \ldots \vee Y_{p-1},
\end{equation}
where the bottom cell of $Y_i$ is in degree $2i+1$.
\begin{rem}\label{rem:Yi}
Each space $Y_i$ inherits a co-$H$-space structure from $\Sigma\CPi_{(p)}$
via the inclusion and projection maps. However, this co-$H$-structure is
not necessarily co-associative and neither are the inclusion and projection
maps necessarily co-$H$-maps. Furthermore, each $Y_i$ is minimal atomic in
the sense of~\cite{AB-PM}, and cannot be equivalent to a suspension except
in the case $i = p-1$ when it does desuspend~\cite{BG-NR}.
\end{rem}

Recent work of Selick, Theriault and Wu~\cite{STW} establishes a
Hilton-Milnor like splitting of loops on such a wedge of co-$H$-spaces.
In our case, the splitting is of the familiar form
\[
\Omega \Sigma \CPi_{(p)} \simeq \Omega(Y_1 \vee \ldots \vee Y_{p-1})
                         \simeq \prod_{\alpha}\!{}^{'}\Omega N_\alpha,
\]
where $\alpha$ runs over all basic products formed on the alphabet
$\{1,\ldots, p-1\}$ and $\prod\!{}^{'}$ denotes the weak product. In
cohomology, \ie, in $\qsym(\Z_{(p)})$, this splitting gives rise to
a splitting of algebras,
\[
H^*(\Omega \Sigma \CPi_{(p)}) \iso \bigotimes_{\alpha} H^*(\Omega N_\alpha).
\]

If $\nu_i(\alpha)$ denotes the number of occurrences of the letter~$i$
in the word $\alpha$, then
\[
\Sigma^{\nu_1(\alpha) + \ldots + \nu_{p-1}(\alpha) - 1}N_\alpha
\simeq
Y_1^{(\nu_1(\alpha))} \wedge \ldots \wedge Y_{p-1}^{(\nu_{p-1}(\alpha))},
\]
where $X^{(n)}$ denotes the $n$th smash power of $X$. The homology of
the space $Y_i$ starts with a generator in degree $2i +1$, thus the
smash power
$Y_1^{(\nu_1(\alpha))} \wedge \ldots \wedge Y_{p-1}^{(\nu_{p-1}(\alpha))}$
has homology starting in degree $\sum (2i+1)\nu_i(\alpha)$. Thus
$\Omega N_\alpha$ has bottom degree
$2\nu_1(\alpha) + \ldots + 2(p-1)\nu_{p-1}(\alpha)$ which is the cohomological
degree of the element $\alpha$ in $\qsym(\Z_{(p)})$.

The map $\CPi \lra BU$ corresponds to the $K$-theory orientation of $\CPi$
in $ku^2(\CPi)$. As $BU \simeq \Omega SU$, there is an adjoint map
$\Sigma \CPi \lra SU$. Since the Adams splitting is compatible with the
infinite loop space structure on $BU$, the delooping of $BU_{(p)}$,
$\B BU_{(p)} \simeq SU_{(p)}$, splits into delooped pieces
$\B W_1 \times \ldots \times \B W_{p-1}$. For a fixed~$j$ in the range
$1 \leq j \leq p-1$, the homology generators
$\beta_i \in H_{2i}(\CPi;\Z_{(p)})$ with $i \equiv j \mod (p-1)$ stem
from $H_*(Y_j;\Z_{(p)})$. The orientation maps $\beta_i$ to the $i$th
generator $b_i$ in $H_*(BU;\Z_{(p)})$ and this lives on the corresponding
Adams summand.

\section{Steenrod operations on $\qsym(\F_p)$}
\label{sec:steenrod}

Working with $\F_p$-coefficients we can ask how the Steenrod operations
connect the generators in $\qsym(\F_p)$. We first state an easy result
about $p$th powers, which is stated without proof in~\cite[(7.17)]{H01}.
\begin{lem}\label{lem:SteenrodOp}
The $p$th power of an element $[a_1,\ldots, a_r] \in \qsym(\F_p)$
is equal to $[pa_1,\ldots, pa_r]$, \ie,
$\mathcal{P}^{a_1 + \ldots + a_r}([a_1,\ldots, a_r]) = [pa_1,\ldots, pa_r]$.
\end{lem}
\begin{proof}
Recall that $[a_1,\ldots, a_r]$ is dual to the generator
$Z_{a_1}\cdot \ldots \cdot Z_{a_r}$ with respect to the monomial basis
in $\nsym$. Thus we have to determine the element which corresponds to
$((Z_{a_1}\cdot \ldots \cdot Z_{a_r})^*)^p$.

Setting $Z(t) = \sum_i Z_i t^i$, we have $\Delta(Z_m) = \sum_i Z_i \otimes Z_{m-i}$,
so $Z(t)$ is group-like. We obtain that the $p$-fold iterate $\Delta^p$ of
the coproduct on a monomial $Z_{j_1}\cdot \ldots \cdot Z_{j_n}$ is captured
in the series
\[
\Delta^p(Z(t_1)\cdot \ldots \cdot Z(t_n))
               = \Delta^p(Z(t_1)) \cdot \ldots \cdot \Delta^p(Z(t_n))
\]
which can be expressed as the $p$-fold product
\[
(Z(t_1)\otimes \ldots \otimes Z(t_n)) \cdot \ldots
                            \cdot (Z(t_1) \otimes \ldots \otimes Z(t_n))
                          = (Z(t_1)\cdot \ldots \cdot Z(t_n))^{\otimes p}.
\]
Thus
\[
\sum_{(j_1,\ldots,j_n)} \langle [a_1,\ldots, a_r]^p, Z_{j_1} \cdot
\ldots \cdot Z_{j_n} \rangle t_1^{j_1}\cdot \ldots \cdot t_r^{j_n}
                           = t_1^{pj_1} \cdot \ldots \cdot t_r^{pj_r},
\]
therefore $[a_1,\ldots, a_r]^p$ is dual to $Z_{pa_1} \clc Z_{pa_r}$.
\end{proof}

We will determine the Steenrod operations in $\qsym$ by using the
James splitting~\eqref{eqn:splitting} and the isomorphism of
Theorem~\ref{thm:qsym}: there is an isomorphism of modules over the
Steenrod algebra
\[
\tilde{H}^*(\lsc; \F_p) \cong \bigoplus_n \tilde{H}^*((\CPi)^{(n)};\F_p).
\]

Let $x_{(i)} \in H^2(\CPi_{(i)})$ be $c_1(\eta_i)$, where $\CPi_{(i)}$
denotes the $i$-th copy of $\CPi$ in the product $(\CPi)^{\times n}$
and $\eta_i$ is the line bundle induced from $\eta$ over $\CPi$. As
$\beta_n$ is dual to $c_1^n$, the elements
$x_{(1)}^{a_1} \.\ldots \.x_{(n)}^{a_n}$ give an additive basis of
$H^*(\lsc)$ and they correspond to the generators $\alpha =[a_1,\ldots,a_n]$.
Since the James splitting~\eqref{eqn:splitting} is only defined after
one suspension, we cannot read off the multiplicative structure of
$H^*(\lsc)$ immediately.

For two generators $\alpha$ and $\beta $ in $\qsym$ we denote their
concatenation by $\alpha \ast \beta$. Thus for $\alpha$ as above and
$\beta = [b_1,\ldots,b_m]$, we have
\[
\alpha \ast \beta = [a_1,\ldots,a_n,b_1,\ldots,b_m].
\]
The Cartan formula leads to a nice description of the action of the
Steenrod algebra on $\qsym(\F_p)$.

\begin{prop}\label{prop:star}
The following Cartan formula holds for the star product:
\begin{equation}\label{eq:star}
\mathcal{P}^i(\alpha \ast \beta) =
\sum_{k+\ell = i} \mathcal{P}^k(\alpha) \ast \mathcal{P}^\ell(\beta).
\end{equation}
\end{prop}
\begin{proof}
Identifying $\alpha$ and $\beta$ with their corresponding elements in
the cohomology of a suitable smash power of $\CPi$, the Cartan formula
on $H^*({(\CPi)}^{(n+m)}; \F_p)$ becomes
\[
\mathcal{P}^i(x_{(1)}^{a_1} \clc x_{(n)}^{a_n} x_{(n+1)}^{b_1}
                                                  \clc x_{(n+m)}^{b_m})
= \sum_{k+\ell=i} \mathcal{P}^k(x_{(1)}^{a_1} \clc x_{(n)}^{a_n})
                  \mathcal{P}^\ell(x_{(n+1)}^{b_1} \clc x_{(n+m)}^{b_m}).
\]
Identifying $x_{(n+i)}^{b_i}$ with $x_{(i)}^{b_i}$ gives the result.
\end{proof}

As we can write every element $\alpha = [a_1,\ldots,a_n]$ as
$[a_1]\ast \ldots \ast [a_n]$ it suffices to describe the Steenrod
operations on the elements of the form $[n]$ with $n \in \N$. The
following is easy to verify using standard identities for binomial
coefficients mod~$p$ and from Lemma~\ref{lem:SteenrodOp}. We
recall that when $s>r$, $\dbinom{r}{s} = 0$.

\begin{prop}\label{prop:Pk[n]}
If the $p$-adic expansions of $n$ and $k$ are $n_0 + n_1p+ \ldots +n_rp^r$
and $k_0 + k_1p+ \ldots + k_sp^s$ respectively, where $s \leq r$, then
\begin{align*}
\mathcal{P}^k[n] & = \binom{n}{k}[n+k(p-1)] \\
 & = \binom{n_0}{k_0} \ldots \binom{n_s}{k_s} [n+k(p-1)].
\end{align*}
Hence if $p \nmid (k_0-n_0)$ and if $n_i \geq k_i$ for all $1\leq i\leq s$,
then $\mathcal{P}^k[n]$ is a non-zero indecomposable. When $p\nmid(k_0-n_0)$
but $n_i < k_i$ for some $i$, the right hand side is zero, and for $p\mid(k_0-n_0)$
the right hand side is decomposable.
\end{prop}
For example, the power operations $\mathcal{P}^k$ on the elements
$\alpha=[a_1, \ldots, a_n]$ with $1 \leq a_i \leq p-1$ and $k\leq\deg(\alpha)$
yield non-trivial sums of indecomposables.

\section{Witt vectors and the cohomology of $BU$}
\label{sec:Witt&LambdaRings}

The self-dual bicommutative Hopf algebras $H_*(BU) \cong H^*(BU)$ are
closely related to both lambda rings and Witt vectors (see~\cite{H0}).
In particular, there are $p$-local splittings due to Husem\"oller~\cite{Hu},
subsequently refined to take into account the Steenrod actions in~\cite{AB-HW}.
Here is a brief account of this theory over any commutative ring~$R$.

First consider the graded commutative Hopf algebra
\[
H_*(BU;R) = R[b_n \;|\; n \geq 1],
\]
where $b_n \in H_{2n}(BU;R)$ is the image of the canonical generator
of $H_{2n}(\CP^\infty;R)$ as defined in~\cite{Ad-0}, and the coproduct
is determined by
\[
\Delta(b_n) = \sum_{i+j=n} b_i \otimes b_j.
\]
This coproduct makes the formal power series
\[
b(t) := \sum_{n\geq0} b_n t^n\in H_*(BU)[[t]]
\]
grouplike, \ie,
\begin{align*}
\Delta b(t) & = b(t)\otimes b(t), \\
\intertext{or equivalently}
\sum_{n\geq0} \Delta(b_n) t^n
 & = \Bigl(\sum_{n\geq0} b_m t^m\Bigr)\otimes\Bigl(\sum_{n\geq0} b_n t^n\Bigr).
\end{align*}
There is an obvious isomorphism of graded Hopf algebras over $R$,
\[
H_*(BU;R) \iso \sym(R)
\]
under which $b_n \leftrightarrow c_n$.

For each $n\geq1$, there is a cyclic primitive submodule
\[
\Prim(\sym(R))_{2n} = R\{q_n\},
\]
where the generators are defined recursively by $q_1=b_1$ together
with the Newton formula
\begin{equation}\label{eqn:NewtonRec}
q_n=b_1q_{n-1}-b_2q_{n-2}+ \ldots + (-1)^{n}b_{n-1}q_1 + (-1)^{n-1}nb_n.
\end{equation}
The $q_n$ can be defined using generating functions and logarithmic
derivatives as follows.
\begin{equation}\label{eqn:qlogder}
\sum_{n \geq 1} (-1)^{n-1}q_nt^{n-1} = \frac{d}{dt}\log(b(t))
                                         =  \frac{b'(t)}{b(t)}.
\end{equation}
The recursion \eqref{eqn:NewtonRec} then follows via multiplication
by $b(t)$. The fact that the $q_n$ are primitive follows because the
logarithmic derivative maps products to sums.

Now we introduce another family of elements $v_n\in\sym(R)_{2n}$
which are defined by generating functions using an indexing
variable~$t$ through the formula
\begin{equation}\label{eqn:WittVect-GenFunc}
\prod_{k\geq1}(1-v_kt^k) = \sum_{n\geq0} b_n (-t)^n.
\end{equation}
These also satisfy
\begin{equation}\label{eqn:WittVect-Rec}
q_n = \sum_{k\ell = n} k v_k^\ell.
\end{equation}
\begin{thm}\label{thm:WittVect}
The elements $v_n\in\sym(R)_{2n}$ are polynomial generators
for $\sym(R)$,
\[
\sym(R) = R[v_n \;|\; n \geq 1].
\]
The coproduct is given by the recursion
\[
\sum_{k\ell = n} k \Delta(v_k)^\ell =
     \sum_{k\ell = n} k (v_k^\ell\otimes1+1\otimes v_k^\ell).
\]
\end{thm}

When $R$ is a $\Z_{(p)}$-algebra, for each $n$ such that $p \nmid n$,
there are elements $v_{n,r} \in \sym(R)_{2np^r}$ defined recursively
by
\[
q_{np^r} =
p^rv_{n,r} +  p^{r-1}v_{n,r-1}^{p} + \ldots + v_{n,0}^{p^r}.
\]
Then the subalgebra
\[
R[v_{n,r} \;|\; r \geq 0] \subset \sym(R)
\]
is a sub Hopf algebra. The following result was first introduced
into topology in~\cite{Hu}.
\begin{thm}\label{thm:WittVect-plocal}
If $R$ is a $\Z_{(p)}$-algebra, there is a decomposition of Hopf
algebras
\[
\sym(R) = \bigotimes_{p \nmid n} R[v_{n,r} \;|\; r \geq 0].
\]
\end{thm}

Notice that when $R$ is an $\F_p$-algebra, we have
$q_{np^r} = q_n^{p^r}$.

There are Frobenius and Verschiebung Hopf algebra endomorphisms
\[
\mathbf{f}_d,\mathbf{v}_d\:\sym(R) \lra \sym(R)
\]
given by
\[
\mathbf{f}_d(v_n)=v_{nd},
\quad
\mathbf{v}_d(v_n)=
\begin{cases}
dv_{n/d}& \text{if $d\mid n$}, \\
0& \text{otherwise}.
\end{cases}
\]

In the dual $H^*(BU;R)$, we have the universal Chern classes $c_i$
and the primitives $s_i$ and under there is an isomorphism of Hopf
algebras
\[
H_*(BU;R)\xrightarrow{\iso}H^*(BU;R);
\quad
b_i \leftrightarrow c_i, \ q_i \leftrightarrow s_i.
\]
We define the element $w_i$ to be the image of the element $v_i$
under this isomorphism. When localized at a prime~$p$, we define
$w_{n,r}$ to be the element corresponding to $v_{n,r}$. The coproduct
on the $w_i$ is computed using the analogue of
Equation~\eqref{eqn:WittVect-GenFunc},
\begin{equation}\label{eqn:WittVect-GenFunc-w}
\prod_{k\geq1}(1-w_kt^k) = \sum_{n\geq0} c_n (-t)^n
\end{equation}
together with the Cartan formula for the $c_i$.

The $c_i$ can be identified with elementary symmetric functions in
infinitely many variables, say $x_i$, and the coproduct $\psi_\oplus$
on a symmetric function $f(x_i)$ amounts to splitting the variables
into two infinite collections, say $x'_i$, $x''_i$, and expressing
the symmetric function $f(x'_i,x''_j)$ in terms of symmetric functions
of these subsets. There is a second coproduct $\psi_\otimes$ corresponding
to replacing $f(x_i)$ by $f(x'_i+x''_j)$. There is an interpretation of
this structure in terms of symmetric functions. For example, the latter
coproduct gives
\[
\psi_\otimes(s_n) = \sum_{0\leq i\leq n}\binom{n}{i}s_i\otimes s_{n-i}.
\]

Both of these coproducts are induced by topological constructions.
Let us recall that the space $BU$ admits maps that represent the
Whitney sum and tensor products of bundles,
\[
BU \times BU \xrightarrow{\oplus} BU,
\quad
BU \times BU \xrightarrow{\otimes} BU.
\]
Using the Splitting Principle, it is standard that the resulting coproducts
\[
H^*(BU)\lra H^*(BU)\otimes H^*(BU)
\]
induced by these are equal to our two coproducts $\psi_\oplus,\psi_\otimes$.

Let $R$ be again a commutative ring with unit. There are two endofunctors
$\Lambda(-)$ and $W(-)$ on the category of commutative rings, $\mathbf{Rings}$.
Details of this can be found in~\cite[chapter~III]{H0}, in particular in E.2,
although the construction there is of an inhomogeneous version corresponding
to $K^0(BU)$ rather than $H^*(BU)$. We will first describe their values in
the category of sets.
\begin{defn}\label{defn:Lambda-Witt}
Let $\Lambda(R)$ be $1 + tR[[t]]$ with addition given by the multiplication
of power series and let the (big) Witt vectors on $R$, $W(R)$, be the
product $\prod_{i \geq 1} R$.
\end{defn}
Consider the following representable functor $\mathbf{Rings}(\sym,-)$.
The coproducts $\psi_\oplus$ and $\psi_\otimes$ make this into a ring
scheme. Therefore they induce two different commutative multiplications
on $\mathbf{Rings}(\sym, R)$.

Let $\varphi_\Lambda\: \mathbf{Rings}(\sym, R) \lra \Lambda(R)$ be the
bijection that sends a map $f$ from $\sym$ to $R$ to the power series
\[
1 + f(c_1)(-t) + f(c_2) (-t)^2 + \ldots + f(c_n) (-t)^n + \ldots.
\]
The coproduct $\psi_\oplus$ corresponds to the multiplication of power
series which should be thought of as a kind of addition, whereas
$\psi_\otimes$ gives a multiplication. These two operations interact
to make this a functor with values in commutative rings.

For the Witt vectors, $W(R)$, take the bijection
$\varphi_W\: \mathbf{Rings}(\sym, R) \lra W(R)$ that sends $f$
to the sequence $(f(w_i))_{i \geq 1}$ where the $w_i$ are given
as in~\eqref{eqn:WittVect-GenFunc-w}. Therefore, the  two
coproducts $\psi_\oplus$ and $\psi_\otimes$ induce a ring
structure on $W(R)$.
\begin{thm}\label{thm:Lambda-Witt}
The two different $H$-space structures, $BU_\oplus$ and $BU_\otimes$
give rise to two comultiplications on\/ $\sym$ via $\psi_\oplus$
and $\psi_\otimes$ which together induce a ring structure on the
big Witt vector $W(R)$ ring.
This ring structure coincides with the standard one as it is
described for instance in~\cite[III \S 17]{H0}.
\end{thm}
\begin{proof}
The exponential map \cite[(17.2.7)]{H0} sends a Witt vector
$(a_1,a_2,\ldots) \in W(R)$ to the element
$\prod_{i \geq 1}(1- a_it^i) \in \Lambda(R)$.
Formula~\eqref{eqn:WittVect-GenFunc-w} ensures that this isomorphism
of rings sends the sequence of generators $(w_n)_n$ to the product
$\prod_{i \geq 1}(1 + (-1)^i c_it^i)$. That the ring structure agrees
for $\Lambda(R)$ follows directly from the definition.
\end{proof}

\section{Quasi-Witt vectors}
\label{sec:q-Witt}

In the case of symmetric functions, we interpreted the addition and
multiplication of Witt vectors as coming from the two $H$-space
structures on $BU$. On $\lsc$ we have the ordinary $H$-space structure,
$\mu_\oplus \: \lsc \times \lsc \lra \lsc$, coming from loop addition.
In addition to this, we consider the following construction. Essentially
the same construction already appears in~\cite{KM}, the difference arises
from various choices concerning the Hopf construction.

Let $X$ be an $H$-space. For two loops $f, g \: \mathbb{S}^1 \lra \Sigma X$,
their product $f \diamond g$ is the composition
\[
\xymatrix{{\mathbb{S}^1} \ar[rrdd]_{f \diamond g} \ar[r]^(0.4)g
 & {\mathbb{S}^1 \wedge X} \ar[r]^(0.4){f \wedge \id}
 & {\mathbb{S}^1 \wedge X \wedge X} \ar[d]^{\simeq} \\
{} & {} & {X \ast X} \ar[d]^{\mathcal{H}(\nu_X)} \\
{} & {} & {\Sigma X}
}
\]
in which the last map is the Hopf construction $\mathcal{H}$ on the $H$-space
multiplication $\nu_X$ on $X$ from the reduced join of two copies of $X$,
see~\cite[XI,~\S 4]{Wh}.

\begin{lem} \label{lem:diamond}
The $\diamond$-product defines a map
$\diamond \: \Omega \Sigma X \wedge \Omega \Sigma X \lra \Omega \Sigma X$.
It satisfies the following left distributivity law for loops
$f,g,h \in \Omega \Sigma X$:
\begin{equation}\label{eqn:diamond-leftdist}
f \diamond (\mu_\oplus(g, h)) \simeq \mu_\oplus(f \diamond g,f \diamond h).
\end{equation}
If $i\:X\lra\Omega \Sigma X$ is the natural map, then for any $x\in X$
the following right distributivity law holds:
\begin{equation}\label{eqn:diamond-rightdist}
\mu_\oplus(g,h)\diamond i(x) \simeq \mu_\oplus(g \diamond i(x),h \diamond i(x)).
\end{equation}
A map of $H$-spaces $\ell\: X \lra Y$ is compatible with the
$\diamond$-operations in that it satisfies
\[
\Omega \Sigma (\ell) \diamond \Omega \Sigma (\ell)
               \simeq \Omega \Sigma (\ell) \circ \diamond.
\]
\end{lem}
\begin{proof}
It is obvious that the $\diamond$-product with the constant loop on
either side gives the constant loop again. For the distributivity
law we consider the following diagram
\[
\xymatrix{
{\mathbb{S}^1} \ar[r]^{\mu_\oplus(g,h)} \ar[ddd]_{\mathrm{pinch}}
& {\mathbb{S}^1 \wedge X} \ar[r]^{f \wedge \id} &
{\mathbb{S}^1 \wedge X \wedge X} \ar[dr] & {}\\
{} &{} &{} &{\Sigma X}\\
{}& {} & {} & {\Sigma X \vee \Sigma X}\ar[u]^{\nabla}\\
{\mathbb{S}^1 \vee \mathbb{S}^1} \ar[r]^(0.4){g \vee h} &
{\mathbb{S}^1 \wedge X \vee \mathbb{S}^1 \wedge X} \ar[r]^(0.4)
{f \wedge \id}_(0.4){\vee f \wedge \id} \ar[uuu]^{\nabla}
& {\mathbb{S}^1 \wedge X \wedge X \vee \mathbb{S}^1 \wedge X \wedge X}
\ar[uuu]^{\nabla} \ar[ur] & {}
}
\]
in which all squares commute and the upper composition corresponds
to $f\diamond (\mu_\oplus)$ whereas the lower composition is
$\mu_\oplus(f \diamond g,f \diamond h)$. Hence the diagram commutes
and left distributivity is verified.

A simpler argument gives the right distributivity for loops of the
form $i(x)$. Naturality of the diamond product with respect to maps
of $H$-spaces follows because the Hopf-construction $\mathcal{H}$ on
an $H$-space multiplication $\nu$ satisfies
\[
\mathcal{H}(\nu_Y) \circ (\ell \ast \ell)
                        \sim \Sigma(\ell) \circ \mathcal{H}(\nu_X)
\]
for $\ell$ an $H$-map from $X$ to $Y$~\cite[XI.4]{Wh}.
\end{proof}

For the following, we recall the notion of a \emph{near-ring},
see~\cite{Pilz-Nearrings} for example. Briefly, a near-ring
is a group equipped with a second product which left or right
distributes over the group operation and is associative. Various
extra conditions are sometimes imposed such as the existence
of a multiplicative unit. Dropping the associativity requirement
leads to the notion of a \emph{quasi-near-ring} which we refrain
from using!
\begin{rem}\label{rem:near-ring}
The structure $(\Omega \Sigma X, \mu_\oplus, \diamond)$ for an
associative $H$-space $X$  from Lemma~\ref{lem:diamond} might
be called a left (non-unital, non-associative) near-ring in the
homotopy category and we abbreviate that to \emph{homotopy near-ring}.
By functoriality, $H_*(\lsc) \cong \nsym$ inherits such a structure
in the category of graded $\Z$-coalgebras. Dually, $H^*(\lsc)\cong\qsym$
is a co-near-ring in the category of graded commutative rings.
\end{rem}

\begin{defn}\label{def:q-Witt}
Let $R$ be a commutative ring with unit. We call
\[
QW(R):= \mathbf{Rings}(\qsym, R)
\]
the near-ring of \emph{quasi-Witt vectors on $R$}.
\end{defn}

Note that these quasi-Witt vectors, $QW(R)$ do not agree with Hazewinkel's
noncommutative Witt vectors $M(R)$ of~\cite[p.75]{H1}. Our coproducts
are induced by maps on spaces, hence they are homogeneous. Hazewinkel's
formula~\cite[\S~12]{H1} gives a non-homogeneous coproduct. They might
agree if we used a $K$-theoretic analogue in place of our homological
one.

{}From Selick~\cite[p.84]{S} we can deduce an explicit formula for the
$\diamond$ product on basic loops $i(x)\:t\mapsto[t,x]\in\Sigma X$ for
an associative $H$-space $X$. It will send two of such basic loops $i(x)$
and $i(y)$ to $\mu_\oplus(\widebar{i(y)},i(x) \cdot i(y), \widebar{i(x)})$,
where $\cdot$ is the $H$-product on $X$ and $\widebar{(\ph{-})}$ denotes
loop reversal. This allows us to calculate $Z_i \diamond Z_j$ in $H_*(\lsc)$.
We note that by construction the diamond product is really defined on
reduced homology. When $x$ is a positive degree element of $H_*(\lsc)$,
we have
\begin{equation} \label{eqn:diamond-1}
1 \diamond x = 0 = x \diamond 1.
\end{equation}

Recall that $Z(t) = \sum_{i\geq0} Z_it^i$.
\begin{prop}\label{prop:diamond-CPi}
The diamond product of the two generators $Z_i,Z_j$
is
\begin{equation}\label{eqn:diamond-CPi}
Z_i \diamond Z_j =
\sum_{r=0}^i\sum_{s=0}^j\binom{r+s}{r}\chi(Z_{j-s})Z_{r+s}\chi(Z_{i-r}).
\end{equation}
In terms of generating functions this is equivalent to
\begin{equation} \label{eqn:diamond-CPigen}
Z(s) \diamond Z(t) = Z(t)^{-1}Z(s+t)Z(s)^{-1}.
\end{equation}
\end{prop}
\begin{proof}
Note that in order to apply Selick's formula, we have to use the diagonal
\[
\lsc \times \lsc \lra (\lsc)^4.
\]
On homology this corresponds to taking the coproduct of $Z_i$ and $Z_j$,
thus $Z_i \otimes Z_j$ maps to
\[
\sum_{r=0}^i\sum_{s=0}^j Z_r \otimes Z_{i-r} \otimes Z_s \otimes Z_{j-s}.
\]
We have to switch factors and apply the $\CPi$-$H$-multiplication to
$Z_r \otimes Z_s$ to give the term $\dbinom{r+s}{r}Z_{r+s}$ because
the cohomology of $\CPi$ is primitively generated. We now apply the
loop-inversion antipode in $H_*(\lsc)$ to the remaining factors, then
finally, we use loop multiplications to obtain the result.
\end{proof}

We can calculate products of the form $u\diamond(v_1\ldots v_n)$ with
$u,v_1,\ldots,v_n$ all of positive degree by using left
distributivity~\eqref{eqn:diamond-leftdist} and~\eqref{eqn:diamond-1}.
Here we denote the loop product by juxtaposing homology elements, \ie,
$xy=x\.y$. For example, if the coproduct on $u$ is
\[
\Delta(u)= u\otimes1 + 1\otimes u + \sum_r u'_r\otimes u''_r,
\]
so that $\sum_r u'_r\otimes u''_r$ is the reduced coproduct of~$u$,
then we have
\begin{equation}\label{eqn:diamond-distrib}
u\diamond(v_1v_2) = \sum_r(u'_r\diamond v_1)(u''_r \diamond v_2).
\end{equation}
In particular, for positive degree elements $x,y$, this gives
\begin{equation}\label{eqn:diamond-distribZileft}
Z_i \diamond (xy) = \sum_{r=1}^{i-1} (Z_r \diamond x)(Z_{i-r} \diamond y)
\end{equation}
for any loop product $xy$ of positive degree elements in $H_*(\lsc)$.
Similarly, right distributivity~\eqref{eqn:diamond-rightdist} gives
\begin{equation}\label{eqn:diamond-distribZiright}
(xy) \diamond Z_i = \sum_{r=1}^{i-1} (x \diamond Z_r)(y \diamond Z_{i-r}).
\end{equation}
In particular,
\begin{align}
Z_1 \diamond (xy)  &= 0,
\label{eqn:diamond-Z1rightannih}\\
(xy) \diamond Z_1  &= 0.
\label{eqn:diamond-Z1leftannih}
\end{align}
The first of these is a special case of the more general
\begin{lem}\label{lem:diamond-primdecomp}
If~$u$ is primitive, then for all positive degree elements $x,y$,
\[
u\diamond(xy) = 0,
\]
\ie, the left diamond product with a primitive annihilates decomposables.
\end{lem}
\begin{proof}
This follows immediately from Equations~\eqref{eqn:diamond-1}
and~\eqref{eqn:diamond-distrib}.
\end{proof}

In the non-commutative algebra $\nsym$ there are two natural families
of primitives analogous to the family $(q_n)_{n \in \N}$
from~\eqref{eqn:NewtonRec}. We define
\begin{align}\label{eqn:leftlogder}
\sum_{n \geq 1} (-1)^{n-1} Q_n t^{n-1} & = Z(t)^{-1}Z'(t), \\
\sum_{n \geq 1} (-1)^{n-1} Q'_n t^{n-1} & = Z'(t)Z(t)^{-1},
\label{eqn:rightlogder}
\end{align}
where $Z'(t)$ is the derivative of the series $Z(t)$. These two generating
functions are related through conjugation by the $Z(t)$-series. Multiplication
by $Z(t)$ on the left (resp. right) gives the recursion formulae
\begin{align}\label{eqn:leftNewtonRec}
Q_n & =
Z_1Q_{n-1}-Z_2Q_{n-2}+ \ldots + (-1)^{n}Z_{n-1}Q_1 + (-1)^{n-1}nZ_n, \\
Q'_n & =
Q'_{n-1}Z_1 - Q'_{n-2}Z_2 + \ldots + (-1)^{n}Q'_1Z_{n-1} + (-1)^{n-1}nZ_n.
\label{eqn:rightNewtonRec}
\end{align}
These families agree with Hazewinkel's~\cite[pp.328--329]{H2} up to
sign. As the antipode on the $Z_i$ is given by the generating function
\[
\sum_{i\geq0}\chi(Z_i)t^i = Z(t)^{-1},
\]
we also obtain the formulae
\begin{align}\label{eqn:leftQnchi}
Q_n & = (-1)^{n-1}\sum_{j=1}^n j\chi(Z_{n-j})Z_j, \\
Q'_n & =
(-1)^{n-1}\sum_{j=1}^n jZ_j\chi(Z_{n-j}).
\label{eqn:rightQnchi}
\end{align}

\begin{prop}\label{prop:Qn-diamond}
The primitives $Q_n$ and $Q'_n$ can be expressed as
\[
Q_n = (-1)^{n-1} Z_1 \diamond Z_{n-1},
\quad
Q'_n =  (-1)^{n-1} Z_{n-1} \diamond Z_{1}.
\]
\end{prop}
\begin{proof}
Using~\eqref{eqn:diamond-CPigen} and taking the derivative
with respect to $s$ at $s=0$ we obtain
\[
Z_1 \diamond Z(t) = Z(t)^{-1}Z'(t) - Z_1,
\]
which gives
\[
Z_1 \diamond Z_{n-1} = (-1)^{n-1} Q_n.
\]
Similarly, using the derivative with respect to $t$ at $t=0$
gives
\[
Z(s) \diamond Z_1 = Z'(s)Z(s)^{-1} - Z_1
\]
and we get the corresponding formula for $Q'_n$.
\end{proof}

We can calculate the recursively defined `power'
\[
Z_1^{\diamond n}=Z_1\diamond(Z_1^{\diamond (n-1)}),
\]
which is actually a spherical element since $Z_1$ is spherical.
An easy induction on~$n$ yields
\begin{equation}\label{eqn:Z1^diamondn}
Z_1^{\diamond n} = (-1)^{n-1}(n-1)!\,Q_n.
\end{equation}
Related formulae given in terms of~\eqref{eqn:leftQnchi} occur
in~\cite[corollary~5.2]{KM} and further spherical primitives are
also determined there. We could also replace the above recursion
by the one involving $(Z_1^{\diamond (n-1)})\diamond Z_1$ which
yields $(-1)^{n-1}(n-1)!\,Q'_n$ in place of $(-1)^{n-1}(n-1)!\,Q_n$.

Now we give a modification of the above construction. For us, a
spectrum will mean a collection of spaces $E=\{E_n\}$ and suitably
related homeomorphisms $\sigma_{n,1}\: E_n \lra \Omega E_{n+1}$ with
adjoint evaluation maps $\tilde{\sigma}_{n,1}\: \Sigma E_n \lra E_{n+1}$.
We obtain maps $\sigma_{n,m}\: E_n \lra \Omega^{m}E_{n+m}$ as iterations
of the maps $\sigma_{n,1}$ and we denote the inverse of $\sigma_{n,m}$
by $\phi_{n+m,n}$.

Such structure is known to be present for many important ring spectra
such as those of complex $K$-theory which is our main concern, it is
always present when $E$ is an $E_\infty$ ring spectrum. We also write
$E'=\{E'_n\}$ for the $0$-connected cover of $E$, so $E'_0$ is the
connected component of the basepoint in $E_0$.

Let $U_1(E)$ denote the component of $E_0$ corresponding to~$1$ in the
ring $\pi_0(E)=\pi_0(E_0)$. Then $U_1(E)$ admits the structure of an
$H$-space under the restriction of $\mu_{0,0} \: E_0 \wedge E_0 \lra E_0$,
and we denote this multiplication by $\mu$.

There is a homotopy equivalence $\lambda\:U_1(E)\lra E'_0$ which satisfies
\begin{equation}\label{eqn:lambda}
\lambda \circ \mu \simeq
(\mu\circ\lambda\times\lambda)*(\lambda\circ\mathrm{pr}_1)*(\lambda\circ\mathrm{pr}_2),
\end{equation}
where $*$ denotes the sum in $E_0$ (loop product) and
$\mathrm{pr}_1,\mathrm{pr}_2$ are the two projections from
$U_1(E) \times U_1(E)$. This is defined by shifting component by adding
an element in the $(-1)$-component of $E_0$. In terms of the cohomology
ring $E^0(X)$ for a connected space $X$, this says that the associated
natural transformation
\[
\bar\lambda\: U_1(E)^0(X) = 1+ \tilde E^0(X) \lra \tilde E^0(X)
\]
satisfies
\[
\bar\lambda(xy) = \bar\lambda(x)\bar\lambda(y)+\bar\lambda(x)+\bar\lambda(y).
\]
We can endow $E'_0$ with the $H$-structure under which $\lambda$ is an
$H$-space equivalence. Then we obtain an induced equivalence of homotopy
near-rings
\[
\Omega\Sigma(\lambda)\: \Omega\Sigma U_1(E) \lra \lse'_0.
\]

The evaluation map $\tilde{\sigma}_{0,1} \: \Sigma E'_0 \lra E'_1$ induces
a loop map
\[
\theta_0 := \phi_{1,0} \circ \Omega\tilde{\sigma}_{0,1}\: \lse'_0 \lra E'_0,
\]
where $E'_0$ forms a homotopy ring object with loop addition and
$\mu' \: E'_0 \wedge E'_0 \lra E'_0$, the restriction of the multiplication
to the $0$-component of $E_0$.
\begin{thm}\label{thm:diamond-E}
The map $\theta_0\: \lse'_0 \lra E'_0$ induces a homotopy commutative
diagram of loop spaces
\[
\xymatrix{
{\Omega\Sigma U_1(E) \wedge \Omega\Sigma U_1(E)}
\ar[r]^(0.6){\diamond} \ar[d]_{\Omega\Sigma\lambda \wedge \Omega\Sigma\lambda}
& {\Omega\Sigma U_1(E)} \ar[d]^{\Omega\Sigma\lambda} \\
{\lse'_0 \wedge \lse'_0} \ar[r]^(0.6){\diamond}
                \ar[d]_{\theta_0 \wedge \theta_0}
 & {\lse'_0}\ar[d]^{\theta_0} \\
{E'_0 \wedge E'_0} \ar[r]^{\mu'} & {E'_0}
}
\]
and hence $\theta_0$ and the composition $\theta_0\circ \Omega\Sigma(\lambda)$
induce maps of homotopy near-rings.
\end{thm}
\begin{proof}
We already know that the top square commutes. For the lower square we
consider the following diagram.
\[
\xymatrix@C=0.5cm{ {\mathbb{S}^2 \cong \mathbb{S}^1 \wedge  \mathbb{S}^1}
\ar[rr]^{\alpha \wedge \beta\quad} \ar[d]^{\id \wedge \beta}&&
{\mathbb{S}^1\wedge E'_0 \wedge  \mathbb{S}^1 \wedge
E'_0}\ar[r]^(0.6){\tilde{\sigma}_{0,1} \wedge \tilde{\sigma}_{0,1}}
   & {E'_1 \wedge E'_1} \ar[r]^{\mu_{1,1}} & {E'_2} \\
{\mathbb{S}^1 \wedge \mathbb{S}^1 \wedge E'_0}
\ar[rr]^{\alpha \wedge \id \wedge \id \ph{abc}}
\ar@/_1pc/[dd]_{\tau \wedge \id}
\ar@/^1pc/[dd]^{(-1) \wedge \id \wedge \id}
 &&{\mathbb{S}^1\wedge E'_0 \wedge \mathbb{S}^1 \wedge E'_0}
  \ar[u]^{=}&{}&{\mathbb{S}^1 \wedge \mathbb{S}^1 \wedge E'_0}
  \ar[u]_{\tilde{\sigma}_{0,2}} \\
{}&&{}&{\mathbb{S}^1 \wedge \mathbb{S}^1 \wedge E'_0 \wedge E'_0}
\ar[ul]^{\id \wedge \tau \wedge \id}
\ar[ur]^{\id \wedge \id \wedge \mu_{0,0}}&{} \\
{\mathbb{S}^1 \wedge \mathbb{S}^1 \wedge E'_0} \ar[rr]^{\id \wedge \alpha \wedge \id \ph{abc}}
&&{\mathbb{S}^1 \wedge \mathbb{S}^1 \wedge E'_0 \wedge E'_0}
\ar@/^1pc/[ur]^(0.3){(-1) \wedge \id \wedge \id \wedge \id}
\ar@/_1pc/[ur]_(0.7){\tau \wedge \id \wedge \id}
&{}&{} }
\]
As we have that
\[
(\id \wedge \tau \wedge \id) \circ (\tau \wedge \id \wedge \id)
\circ (\id \wedge \alpha \wedge \id) \circ (\tau \wedge \id)
= \alpha \wedge \id \wedge \id,
\]
the lower left
outer
pentagon commutes. The maps $\tau \wedge \id \wedge \id$
and $(-1) \wedge \id \wedge \id \wedge \id$ are homotopic thus the lower
left
inner
pentagon commutes up to homotopy. The top horizontal map is
the smash product of $\alpha$ and $\beta$ followed by evaluation and
multiplication. The upper left square commutes and the upper right pentagon
commutes up to homotopy. Therefore this composite is homotopic to
\[
\tilde{\sigma}_{0,2} \circ (\id \wedge \id \wedge \mu_{0,0}) \circ
((-1) \wedge \id \wedge \id \wedge \id) \circ (\id \wedge \alpha \wedge \id)
\circ ((-1) \wedge \id \wedge \id) \circ (\id \wedge \beta)
\]
and this composite is $\tilde{\sigma}_{0,2} \circ \Sigma (\alpha \diamond \beta)$.
\end{proof}

\begin{cor}\label{cor:diamond-E}
Let $R$ be a commutative ring for which $H_*(E'_0;R)$ is $R$-flat.
Then the maps $\theta_0$ and $\theta_0 \circ \Omega\Sigma(\lambda)$
induce epimorphisms of near-ring objects in $R$-coalgebras
\[
H_*(\lse'_0;R) \lra H_*(E'_0;R) \lla  H_*(\Omega\Sigma U_1(E);R).
\]
\end{cor}
\begin{proof}
This follows immediately from the fact that the inclusion map
$E'_0 \lra \lse'_0$ is a right homotopy inverse for $\theta_0$.
\end{proof}

Returning to our main example, we take $E=ku$ and $E'=\Sigma^2 ku$,
thus $E_0=BU\times\Z$ and $E'_0=BU$. Notice that we have a map of
$H$-spaces $\nu\: \CPi \lra U_1(KU)=BU_\otimes$ classifying line
bundles. Our main conclusion is contained in
\begin{thm}\label{thm:diamond-ECP-BU}
The map
$\Theta=\theta_0\circ \Omega\Sigma(\lambda\circ\nu)\: \CPi \lra BU$
induces a commutative diagram
\[
\xymatrix{
{\Omega\Sigma \CPi \wedge \Omega\Sigma \CPi}
\ar[r]^(0.6){\diamond} \ar[d]_{\Theta \wedge \Theta}
& {\Omega\Sigma \CPi} \ar[d]^{\Theta} \\
{BU \wedge BU} \ar[r]^(0.6){\otimes} & {BU}
}
\]
and hence $\Theta$ induces a map of homotopy near-rings. Furthermore,
$\Theta$ induces an epimorphism $\Theta_*\:H_*(\lsc) \lra H_*(BU)$
of near-ring objects in $\Z$-coalgebras.
\end{thm}

Of course, $H_*(BU)$ is already a (non-unital) ring object in $\Z$-coalgebras.

\section{The Thom spectrum of $\xi$}\label{sec:ThomSpectrum-xi}

The map $j\:\lsc \lra BU$ is a loop map and hence the Thom spectrum
$M\xi$ of the virtual bundle $\xi$ is an $A_\infty$ ring spectrum and
the natural map $M\xi \lra MU$ is one of $A_\infty$ ring spectra, or
equivalently of $\mathbb{S}$-algebras in the sense of~\cite{EKMM}. Since
the homology $H_*(M\xi)$ is isomorphic as a ring to $H_*(\lsc)$, we see
that $M\xi$ is not even homotopy commutative, let alone an $E_\infty$
ring spectrum or equivalently a commutative $\mathbb{S}$-algebra.

Let $\mathcal{A}_p^*$ be the mod~$p$ Steenrod algebra and let
$\mathcal{A}_p^*/(Q^0)$ be the quotient by the two-sided ideal generated
by the Bockstein $Q^0$.

We also have in all cases, $H^*(BP;\F_p) \iso \mathcal{A}_p^*/(Q^0)$ as
$\mathcal{A}_p^*$-modules, where $BP$ is the $p$-primary Brown-Peterson
spectrum~\cite[II \S 16]{Ad-1}. More generally, any $\mathcal{A}_p^*$-module
on which $Q^0$ acts trivially can be viewed as a $\mathcal{A}_p^*/(Q^0)$-module.
In particular, $H^*(MU;\F_p)$ becomes a free $\mathcal{A}_p^*/(Q^0)$-module in
this way, and therefore $MU_{(p)}$ splits as a wedge of suspensions of $BP$'s.
The proof of this uses a result of Milnor and Moore~\cite[theorem~4.4]{M&M}
which only requires the existence of an associative coalgebra structure.
Although $M\xi$ is not a commutative ring spectrum we may still apply such an
argument to its cohomology.
\begin{prop}\label{prop:M&M-BP}
Let $p$ be a prime. Then $H^*(M\xi;\F_p)$ is a free
$\mathcal{A}_p^*/(Q^0)$-module and hence $M\xi_{(p)}$ is a wedge of $BP$'s.
\end{prop}
\begin{proof}
The standard Milnor-Moore argument works since $M\xi$ is a ring spectrum
and the map $Mj\:M\xi \lra MU$ induced from $j$ gives rise to a
monomorphism $Mj\: H^*(MU;\F_p) \lra H^*(M\xi;\F_p)$. From this we can
deduce that the copy of $\mathcal{A}_p^*/(Q^0)$ on the stable Thom
class of $H^*(MU;\F_p)$ maps isomorphically to one in $H^*(M\xi;\F_p)$.
Obstruction theory now leads to the existence of an equivalence
$M\xi_{(p)} \lra \bigvee_\lambda \Sigma^{2k_\lambda}BP$.
\end{proof}

Of course, as described in~\cite{Ad-1}, Quillen gave a multiplicative
idempotent on $MU_{(p)}$ leading to a precise description of such a
splitting into $BP$'s. We expect that $M\xi$ has an adequate analogue
of the universality for complex orientations possessed by $MU$ and
required for Quillen's approach. We will investigate this in future
work. We also remark that although the spectrum $M\xi$ appears to be
related to that studied in~\cite{A&B}, the latter is not a ring spectrum.

Instead, we note some algebraic facts that are useful here. Fixing a
prime~$p$, we set $\mathcal{A}^*=\mathcal{A}_p^*$ and write $\mathcal{A}_*$
for the dual Steenrod algebra. The dual of $\mathcal{A}_p^*/(Q^0)$ is
the polynomial sub-Hopf algebra
$\F_p[\zeta_i \;|\; i\geq1] \subset \mathcal{A}_*$ generated by the
conjugates $\zeta_i$ of the Milnor generators $\xi_i$ and there is
an isomorphism of $\mathcal{A}_*$-comodule algebras
\[
H_*(BP;\F_p) \iso \F_p[\zeta_i \;|\; i\geq1].
\]
There is also an isomorphism of $\mathcal{A}_*$-comodules
\[
H_*(M\xi;\F_p) \iso
        \bigoplus_\lambda \Sigma^{2k_\lambda}\F_p[\zeta_i\;|\;i\geq1].
\]

For any commutative ring $R$, under the Thom isomorphism
\[
H_*(\lsc;R) \iso H_*(M\xi;R),
\]
the generator $Z_i$ corresponds to an element $z_i \in H_{2i}(M\xi;R)$
(we set $z_0=1$). Thomifying the map $i\: \CPi \lra \lsc$, we obtain
a map $Mi\: \Sigma^\infty MU(1) \lra \Sigma^2M\xi$ and it is easy to
see that
\begin{equation}\label{eqn:beta-z}
Mi_*\beta_{i+1}=z_i.
\end{equation}
Working in mod~$p$ homology and using the standard formula for the
(left) coaction $\psi$~\cite{Ad-1}, we obtain
\[
\sum_{i\geq0}\psi(z_{i})t^{i+1} = \sum_{j\geq0}\xi(t)^{j+1} \otimes z_{j},
\]
where
\[
\xi(t) = \sum_{k\geq0} \xi_k t^{p^k}.
\]
We now introduce elements $w_i \in H_{2i}(M\xi;R)$ by requiring that
they satisfy the functional equation
\begin{equation}\label{eqn:FuncEqn}
\sum_{i\geq0} z_i(\sum_{j\geq0} w_j t^{j+1})^{i+1} = t,
\end{equation}
where we treat $t$ as a variable that commutes with everything in sight.
This amounts to an infinite sequence of equations of the general form
\begin{equation}\label{eqn:FuncEqn-2}
w_n + \ldots + z_n = 0 \quad (n \geq 1),
\end{equation}
where the intermediate terms involve the elements $z_i,w_i$ with $i=1,\ldots,n-1$
and $w_0 = 1$. Hence we can recursively solve this system of equations for
the $w_n$ and the solution is obviously unique. Mapping into $H_*(MU;R)$ we
obtain the usual generators $m_i \in H_{2i}(MU;R)$ described in~\cite{Ad-1}
and these can be expressed explicitly using the Lagrange inversion formula.
But it is not obvious in what sense such a formula exists within the
non-commutative context where we are working.

Now we focus on the case of $H_*(M\xi;\F_p)$ and investigate the
$\mathcal{A}_*$-coaction $\psi$ on the elements $w_i$. We will make heavy
use of the fact that $\psi$ is multiplicative, \ie, for homogeneous elements
$u \in H_{2r}(M\xi;\F_p)$ and $v \in H_{2s}(M\xi;\F_p)$ with
\[
\psi(u) = \sum_{i}x_i \otimes u_i \in \mathcal{A}_*\otimes H_*(M\xi;\F_p),
\quad
\psi(v) = \sum_{j}y_j\otimes v_j \in \mathcal{A}_*\otimes H_*(M\xi;\F_p),
\]
where $u_i,v_j \in H_*(M\xi;\F_p)$ and $x_i,y_j \in \mathcal{A}_*$,
we have
\[
\psi(uv) = \sum_{i,j}x_iy_j \otimes u_iv_j \in \mathcal{A}_*\otimes H_*(M\xi;\F_p).
\]
Let us write $\tilde w_i = \psi(w_i)$ and
\[
\tilde w(t) = \sum_{i\geq0} \tilde w_i t^{i+1}
                         \in (\mathcal{A}_*\otimes H_*(M\xi;\F_p))[[t]].
\]
Then from~\eqref{eqn:FuncEqn} we have
\begin{equation}\label{eqn:FuncEqn-3}
\sum_{i\geq0} (1 \otimes z_i)[(\xi \otimes 1)(\tilde w(t))]^{i+1} = t,
\end{equation}
where $t$ and all elements of $\mathcal{A}_*\otimes1$ commute with
everything. Of course,
\[
(\xi \otimes 1)(\tilde w(t)) =
                    \sum_{r\geq0}(\xi_r \otimes 1)(\tilde w(t))^{p^r}.
\]
Now unravelling the coefficients we find that the $\tilde w_i$ satisfy
a sequence of equations of same general form as~\eqref{eqn:FuncEqn-2}:
\[
\tilde w_n + \ldots + 1 \otimes z_n = 0 \quad(n \geq 1),
\]
and clearly this has a \emph{unique} solution. Now consider the series
\[
(\zeta \otimes 1)(1 \otimes w)(t) =
            \sum_{r\geq0}(\zeta_r \otimes 1)[(1 \otimes w)(t)]^{p^r}.
\]
Using~\eqref{eqn:FuncEqn}, we obtain
\begin{align*}
\sum_{i\geq0}
 (1 \otimes z_i )[(\xi\otimes 1)(\zeta\otimes 1)(1\otimes w(t))]^{i+1}
 &= \sum_{i\geq0} (1 \otimes z_i )(1 \otimes w(t)^{i+1}) \\
 &= \sum_{i\geq0} 1 \otimes z_i w(t)^{i+1} = t,
\end{align*}
hence we have
\begin{equation}\label{eqn:Coaction-wn}
\sum_{i\geq0} \psi(w_i) t^{i+1} = \sum_{j\geq0} \zeta_j \otimes w(t)^{p^j}.
\end{equation}
This has the same form as the coaction on the $m_i$ in $H_*(MU;\F_p)$
as given in~\cite{Ad-1}, although the non-commutativity means that the
explicit formulae for $\psi(w_n)$ are more complicated. In particular,
the simple formula
\[
\psi(m_{p^r-1})= \sum_{0\leq i\leq r} \zeta_i \otimes  m_{p^{r-i}-1}^{p^{i}}
\]
is replaced by one involving many more terms.

\section{The homology of  $\thh(M\xi)$}
\label{sec:thh}

The following is a description of the result of applying the B\"okstedt
spectral sequence to $M\xi$ and of some things we learned from~\cite{BCS}.
In his influential unpublished preprint~\cite{MB-THH} (see for
instance~\cite[3.1]{MS} for a published account), Marcel B\"okstedt
developed a spectral sequence for calculating the homology of the
$\thh$-spectrum of a strictly associative spectrum.

The spectral sequence for the integral homology of the associative
spectrum $M\xi$ has as its $\mathrm{E}^2$-page the Hochschild homology
of $H_*(M\xi)$, so we have
\[
\mathrm{E}^2_{s,t} = HH_{s,t}(H_*(M\xi)) \Lra H_{s+t}(\thh(M\xi)).
\]
Now $H_*(M\xi)$ is a tensor algebra $TV$, where $V$ is the free
$\Z$-module generated by the elements $z_i$ for $i \geq 1$, thus
$H_*(M\xi)$ is free as an associative algebra and the Hochschild
homology vanishes except for homological degrees zero and one.
Moreover,
\begin{align*}
HH_0(H_*(M\xi)) &= \bigoplus_{m \geq 0}V^{\otimes m}/C_m, \\
HH_1(H_*(M\xi)) &= \bigoplus_{m \geq 1}{(V^{\otimes m})}^{C_m},
\end{align*}
where $V^{\otimes m}/C_m$ and ${(V^{\otimes m})}^{C_m}$ denote
the coinvariants and invariants respectively for the action of
the cyclic group $C_m$. As there are just two non-trivial columns,
this spectral sequence collapses placing the part arising as the
coinvariants in total even degree and the part coming from the
invariants in total odd degree. We note that for each degree~$m$,
a basis of $V^{\otimes m}$ consists of the tensors of length~$m$
in the $z_i$ and this is a permutation basis for the action of $C_m$.
Therefore $V^{\otimes m}$ decomposes as a direct sum of $C_m$-modules
corresponding to the orbits of the action and for each of these it
is straightforward to see that the coinvariants and invariants are
both free of rank one over~$\Z$. Hence $HH_0(H_*(M\xi))$ and
$HH_1(H_*(M\xi))$ are finitely generated free abelian groups in each
bidegree and therefore $H_*(\thh(M\xi))$ consists of such groups in
each degree.

Andrew Blumberg, Ralph Cohen and Christian Schlichtkrull~\cite{BCS}
give a nice description of topological Hochschild homology of Thom
spectra. Let $BF$ denote the classifying space for spherical
fibrations. Given a based map $f\: X \lra \B BF$, which induces a
loop map from the based loop space $\Omega X$ to $BF$, they describe
$\thh(Mf)$ as the Thom spectrum associated to the composite
\[
LX \lra  L\B BF \simeq BF \times \B BF
           \xrightarrow{\id \times \eta} BF \times  BF \lra BF.
\]
Here, $L(-)$ denotes the free loop space on a space. In particular,
the homology of $\thh(Mf)$ is the homology of the free loop space
on $X$ if we actually start with a map to $BSF$.

Let $\mathrm{ad}\:\Sigma \CPi \lra \B BU$ be the adjoint of the
inclusion $j\:\CPi= BU(1) \lra BU$.

In our case, the map to $BF$ factors over $BU$ and thus over $BSF$
and we obtain that $\thh(M\xi)$ is the Thom spectrum associated to
\[
L\Sigma \CPi \xrightarrow{L\mathrm{ad}} L \B BU = L(SU)
         \simeq \Omega SU \times SU = BU \times SU \lra BU
\]
because $\eta$ becomes null homotopic here.

Goodwillie's result~\cite[Theorem~V.1.1]{Go} identifies the
$\mathbb{S}^1$-equivariant homology of the free loop space with
cyclic homology of the chains on the based (Moore) loop space.
The corresponding result for Hochschild homology reads in our
case
\[
H_*(L\Sigma \CPi) \cong HH_*(C_*(\lsc)).
\]
Thus the $\mathrm{E}^2$-term of associated hyperhomology spectral
sequence has the same form as that of the B\"okstedt spectral
sequence. We give a general comparison result for this situation.

Let $X = \Omega Y$ be a loop space  with torsion-free homology and
with a loop map $f \: X \lra BSF$. Assume that $Y$ is path connected
and well-pointed. We can consider the associated B\"okstedt and
hyperhomology spectral sequences for $X$ and its Thom spectrum $Mf$.

\begin{prop}\label{prop:SS-comparison}
There is an isomorphism between the hyperhomology spectral sequence
and the B\"okstedt spectral sequence induced by an isomorphism of
exact couples.
\end{prop}
\begin{proof}
For the hyperhomology spectral sequence we have to consider the
bicomplex with $(p,q)$-term $(\bar{C}_{q}(C_*(\Omega Y)))_p$ where
$\bar{C}_*$ denotes the reduced Hochschild complex. We take the
$s$-filtration of the total complex by considering
\[
F_s(\text{hyper})_n = \bigoplus_{\substack{p+q=n \\ q \leq s}}
(\bar{C}_{q}(C_*(\Omega Y)))_p.
\]
The corresponding exact couple then consists of the terms
\begin{align*}
D(\text{hyper})_{s,t} & = H_{s+t}F_s(\text{hyper})_*,  \\
E(\text{hyper})_{s,t} & = H_{s+t}
(F_s(\text{hyper})_*/F_{s-1}(\text{hyper})_*) =
(\bar{C}_{s}(H_*(\Omega Y)))_t.
\end{align*}

For the B\"okstedt spectral sequence, we use the skeletal
filtration of $\thh(Mf)$. This gives rise to another exact
couple with
\[
D(\text{B\"o})_{s,t} = H_{s+t} (\thh^{[s]}(Mf)),
\]
in which $\thh^{[s]}(Mf)$ denotes the $s$-skeleton of the
realization of the simplicial spectrum $\thh(Mf)$. Since
topological Hochschild homology gives rise to a proper
simplicial spectrum we can apply~\cite[X.2.9]{EKMM} in
order to identify the filtration quotients and get
\[
E(\text{B\"o})_{s,t} \cong H_{s+t}(\Sigma^s (Mf \wedge
\underbrace{\widebar{Mf} \wedge \ldots \wedge \widebar{Mf}}_{s})).
\]
Here, $\widebar{Mf}$ is the cofibre of $\mathbb{S} \lra Mf$.
Therefore we obtain
\[
E(\text{B\"o})_{s,t} \cong (\widebar{C}_{s}(H_*(Mf)))_t.
\]
The Thom isomorphism yields an isomorphism between the two exact
couples and its multiplicativity ensures that the higher differentials
are preserved as well.
\end{proof}

\end{document}